\documentclass[twoside,10.5pt]{article}
\usepackage{mathrsfs}
\usepackage{pifont}
\usepackage{amsmath}
\usepackage{amsthm}
\usepackage{txfonts}
\usepackage{geometry}
\usepackage{latexsym}
\usepackage{amssymb}
\usepackage{graphicx}
\usepackage{geometry}
\usepackage{fancyhdr}
\usepackage{multirow}
\usepackage{xcolor} 
\begin{document}

\newcommand{\li}{\mathop{\mathrm{li}}}
\newcommand{\Gumbel}{\mathop{\mathrm{Gumbel}}}
\newcommand{\Exp}{\mathop{\mathrm{Exp}}}

{\em International Journal of Statistics and Probability}, 3, No.\,2, 18-29, 2014. \ \ ISSN 1927-7032 \ \ E-ISSN 1927-7040
\\
{\em This is arXiv:1401.6959v3 (minor errata fixed; see last page). For the actual journal article, 
 see DOI URL below.}

\begin{center}
{\LARGE{The Distribution of Maximal Prime Gaps in Cram\'er's \\
Probabilistic Model of Primes}}\\[20pt]
\end{center}

\begin{center}
Alexei Kourbatov$^1$
\end{center}
$^1$ JavaScripter.net/math, Redmond, Washington, USA \par
Correspondence: Alexei Kourbatov, 15127 NE 24th Street~\#578, Redmond, WA 98052, USA. 
\\ 
E-mail: akourbatov@gmail.com
\\[15pt]
Received: February 7, 2014  \quad  Accepted: March 11, 2014  \quad Online Published: March 20, 2014 
\\ 
doi:10.5539/ijsp.v3n2p18  \quad \quad  URL: http://dx.doi.org/10.5539/ijsp.v3n2p18 
\\[10pt]



\textbf{Abstract}

In the framework of Cram\'er's probabilistic model of primes, we explore the exact and asymptotic
distributions of maximal prime gaps. We show that the Gumbel extreme value distribution $\exp(-\exp(-x))$
is the limit law for maximal gaps between Cram\'er's random ``primes.'' The result can be derived from a
general theorem about intervals between discrete random events occurring with slowly varying probability  
monotonically decreasing to zero. A straightforward generalization extends the Gumbel limit law 
to maximal gaps between prime constellations in Cram\'er's model.



\medskip
\section{Introduction}

In this paper we apply {\em extreme value theory} to Cram\'er's {\em probabilistic model of primes}. 
But first let us say a few words about two mathematicians who pioneered these important topics.
The Swedish mathematician Harald Cram\'er (1893--1985) made long-lasting contributions in statistics and number theory.
His model of primes continues to serve as a heuristic tool leading 
to new insights in the distribution of primes. Cram\'er wrote:
\begin{quotation}\noindent
In investigations concerning the asymptotic properties of arithmetic
functions, it is often possible to make an interesting heuristic use of
probability arguments. If, e.\,g., we are interested in the distribution of
a given sequence $S$ of integers, we then consider $S$ as a member of 
an infinite class $C$ of sequences, which may be concretely interpreted as
the possible realizations of some game of chance. It is then in many
cases possible to prove that, {\em with a probability} = 1, a certain relation 
$R$ holds in $C$, i.\,e.~that in a definite mathematical sense ``almost all''
sequences of $C$ satisfy $R$. Of course we cannot in general conclude that $R$
holds for the particular sequence $S$, but results suggested in this way
may sometimes afterwards be rigorously proved by other methods.
(Cram\'er, 1936, p.\,25)
%
\end{quotation}
It is difficult to ascertain whether Harald Cram\'er had ever met in person his contemporary, 
the German-American mathematician Emil Julius Gumbel (1891--1966), one of the founders of 
{\em extreme value theory}. 
This branch of statistics is used today for describing phenomena in vastly diverse areas, 
ranging from actuarial science to hydrology to number theory. 
In {\em Statistics of Extremes} (1958) Gumbel observed:
\begin{quotation}\noindent
$\ldots$ many engineers and practical statisticians $\ldots$ are inclined to believe that,
after all, nearly everything should be normal, and whatever turns out not to be so 
can be made normal by a logarithmic transformation. This is neither practical nor true.
(Gumbel, 1958, p.\,345)
%
\end{quotation}
In {\em Les valeurs extr\^emes des distributions statistiques} (1935)
Gumbel showed that extreme values taken from a sequence of 
i.i.d.\,random variables with an {\em exponential} distribution obey
the {\em double exponential} limit law (now known as the Gumbel distribution). 
He reconfirmed the earlier result of Fisher and Tippett (1928) 
that the same limit law also holds for extreme values of i.i.d.\,\,random variables with a normal distribution~--- 
and generalized it to a much wider class of initial distributions, the so called {\em exponential type}. 
To wit:
\begin{quotation}
\noindent
Pour les valeurs extr\^emes $\ldots$ on arrive \`a la distribution doublement 
exponentielle pourvu que la distribution initiale appartienne au type exponentiel 
$\ldots$ Cette theorie est susceptible de nombreuses applications, puisque 
en particulier, la distribution de Gauss et la distribution exponentielle 
appartiennent au type exponentiel. 
(Gumbel, 1935, p.\,154--155) 
\end{quotation}
We will draw upon the work of both Cram\'er and Gumbel: 
In Section 4.2 we show that the limiting distribution of maximal prime gaps 
in Cram\'er's probabilistic model is the Gumbel extreme value distribution. 
In hindsight, the result is not very surprising: we know from computation and distribution fitting 
that the {\em actual} maximal prime gaps are indeed nicely approximated by the Gumbel distribution
(Kourbatov, 2013). 
(What {\em is} somewhat surprising is that Cram\'er or Gumbel have not themselves 
published a similar result long ago, perhaps in the 1930s or 1940s. 
Did they possibly deem it obvious? We might never know.)

The existence of the Gumbel limit law for maxima of the (non-identically distributed)
gaps between Cram\'er's random primes is interesting in view of {\em Mejzler's theorem} 
(see, e.\,g., de Haan and Ferreira, 2006, p.\,201). 
For extremes of non-identically distributed independent random variables, Mejzler's theorem states 
that the limiting distribution (if it exists at all) can be {\em any} distribution with a log-concave cdf.
Thus, prime gaps in Cram\'er's model give us an example of non-identically distributed
random variables whose extremes nevertheless possess a limit law that is
allowed in i.i.d.~situations as well.

\section{Definitions. Notations. Abbreviations}

\begin{tabular}{ll}  
a.s.              & almost sure, almost surely                                                             \\
i.i.d.            & independent and identically distributed                                                \\
cdf               & cumulative distribution function                                                       \\
pdf               & probability density function                                                           \\
$Ex  $            & the {\em expected value} ({\em mathematical expectation}) of the random variable $x$   \\
$\Exp(x;a)$       & the exponential distribution cdf: \ $\Exp(x;a)=1-e^{-x/a}$                             \\
$\Gumbel(x;a,\mu)$& the Gumbel distribution cdf: \ $\Gumbel(x;a,\mu) = 
                    e^{-e^{-{x-\mu\over\vphantom{f}a}}} = 2^{-e^{-{x-M\vphantom{f}\over\vphantom{f}a}}}$   \\
$a$               & the {\em scale parameter} of exponential/Gumbel distributions, as applicable           \\
$\mu$             & the {\em location parameter} ({\em mode}) of the Gumbel distribution                   \\
$M$               & the {\em median} of $\Gumbel(x;a,\mu)$: \ $M\,=\,\mu-a\log\log2\,\approx\,\mu+0.3665a$ \\
RP                & a random ``prime'' in the context of Cram\'er's model (a white ball)                   \\
RC                & a random ``composite'' in the context of Cram\'er's model (a black ball)               \\
$p_k$             & the $k$-th prime;\, $\{p_k\} = \{2,3,5,7,11,\ldots\}$                                  \\
$P_k$             & the $k$-th random ``prime'' (RP) in Cram\'er's model                                   \\
$U_n$             & the $n$-th urn producing RPs with probability ${1\over\log n}$ in Cram\'er's model     \\
$R_n$             & a random variable: the longest uninterrupted run of RCs $\le n$                        \\
$G_n$             & a random variable: the maximal gap between RPs $\le n$; \,\ $G_n:=R_n+1$               \\
$\pi(x)$          & the prime counting function, the total number of primes $p_k\le x$                     \\
$\Pi(x)$          & the RP counting function, the total number of RPs $P_k\le x$                           \\
$\log x$          & the natural logarithm of $x$                                                           \\
$\li x$           & the logarithmic integral of $x$: \
                         $\li x \,= \displaystyle\int_0^x{\negthinspace}{dt\over\log t}
                                \,= \int_2^x{\negthinspace}{dt\over\log t} + 1.04516\ldots$                \\
\end{tabular}

\subsection{Definitions of gaps and runs}

{\it Prime gaps} are distances between two consecutive primes, $p_{k}-p_{k-1}$ (OEIS A001223, Sloane, 2014).
In the context of Cram\'er's model (see next section) ``prime'' gaps refer to distances between consecutive RPs, 
$P_{k}-P_{k-1}$.

{\it Maximal gaps} between primes are usually defined as gaps that are 
strictly greater than all preceding gaps (OEIS A005250, Sloane, 2014).
However, for Cram\'er's model, we will use the term {\em maximal gap} in the statistical sense 
%
%
%
defined below.
Note that Cram\'er's model does not guarantee that there are any ``primes'' $P_k>2$ at all.
To make sure that maximal gaps are defined in all cases, 
we first define the longest run of random ``composites'' (RCs) $\le n$:
$$
R_n ~=~ \mbox{the longest run of consecutive RCs\,} \le n \mbox{ \ (allowing runs of length 0)}.
$$
We now define the {\em maximal gap} between $\mbox{RPs\,}\le n$ simply as
$$
G_n ~=~ R_n + 1.
$$
Compare our definition of $G_n$ to the definition of maximal prime gaps for true primes: 
$$
\mbox{maximal prime gap up to }n\, ~=~ \displaystyle\max_{p_{k}\le n}(p_{k}-p_{k-1}).
$$
Clearly, the same gap/run relation holds (with very rare exceptions) for true primes as well:
$$
\mbox{maximal prime gap} ~=~ 1 ~+~ \mbox{the longest run of composites below } n.
$$
{\em Examples:} the largest prime gap below 100 is the gap of 8~--- between the primes 89 and 97;
this corresponds to a run of 7 consecutive composites: 90 to 96.
Exceptional cases ($gap \ne 1 + run$) occur towards the end of record prime gaps. 
Consider, e.\,g., the prime gap $113\ldots127$ and take $n=125$.
Then the largest prime gap below $n$ is still the gap of 8~--- between 89 and 97~--- 
while the longest run of 12 composites $\le n$ occurs from 114 to 125.

Of course, in the probabilistic model, both $R_n$ and $G_n$ are random variables and,
by definition, $G_n = R_n + 1$ without exceptions.
In Section 4 we will investigate the distribution of maximal gaps $G_n$.

\subsection{Remark: maxima vs.\ records}

Number theorists may use the terms {\it maximal gaps} and {\it record gaps} as synonyms 
(Caldwell, 2010; Nicely, 2013) 
while statisticians make a distinction between {\em maximal} and {\em record} values 
(Arnold {\em et al.}, 1998; Nevzorov, 2001).
Resnick (1973, Theorem 8) shows that the distinction is quite profound: in the i.i.d.\,\,case,
the limit laws for records and maxima {\em cannot be the same}.
Nevertheless, for a wide class of sequences of {\em non}-identically distributed random variables
(the $F^\alpha$\,{\em model}) the same extreme value distribution (e.\,g.~Gumbel distribution) 
can be the limit law for {\em both} records and maxima (Arnold {\em et al.}, 1998, p.\,193).
%
Clearly, Cram\'er's ``prime'' gaps near, say, urn $U_{10}$ are distributed {\em not} identically 
to those near urn $U_{100}$. Theorem 1 (Sect.\,4.2) establishes that the limiting 
distribution of Cram\'er's maximal ``prime'' gaps is indeed the Gumbel distribution. 
Computational evidence suggests that the Gumbel distribution is also the a.s.\,\,limit law 
for record gaps in Cram\'er's model; we will discuss the distribution of records elsewhere.

\section{Cram\'er's probabilistic model of primes}

Cram\'er's probabilistic model of primes is well known~---
and much criticized (Granville, 1995; Pintz, 2007).

\subsection{Setting up the model}

Cram\'er (1936) sets up the model of primes as follows:

\begin{quotation}\noindent{
With respect to the ordinary prime numbers, it is well known that,
roughly speaking, we may say that the chance that a given integer $n$
should be a prime is approximately ${1\over\log n}$. This suggests that by considering
the following series of independent trials we should obtain sequences
of integers presenting a certain analogy with the sequence of ordinary
prime numbers $p_n$.

\noindent
Let $U_1$, $U_2$, $U_3,\ldots$ be an infinite series of urns containing black and
white balls, the chance of drawing a white ball from $U_n$ being ${1\over\log n}$ for
$n > 2$, while the composition of $U_1$ and $U_2$ may be arbitrarily chosen.
We now assume that one ball is drawn from each urn, so that an infinite
series of alternately black and white balls is obtained. If $P_n$ denotes the
number of the urn from which the $n$-th white ball in the series was drawn,
the numbers $P_1$, $P_2,\ldots$ will form an increasing sequence of integers, and
we shall consider the class $C$ of all possible sequences $(P_n)$. Obviously
the sequence $S$ of ordinary prime numbers $(p_n)$ belongs to this class.

\noindent
We shall denote by $\Pi(x)$ the number of those $P_n$ which are $\le x$,
thus forming an analogy to the ordinary notation $\pi(x)$ for the number
of primes $p_n\le x$. (Cram\'er, 1936, pp.\,25--26)
}
\end{quotation}
Cram\'er's model, as stated, is {\em underdetermined}:
the content of urns $U_1$ and $U_2$ is arbitrary.
To compute exact distributions of maximal gaps, we will assume that 

  (i) urn $U_1$ is empty~--- it produces neither ``primes'' nor ``composites'';

 (ii) urn $U_2$ always produces white balls (i.\,e.~the number 2 is certain to be ``prime'').

\subsection{Cram\'er's results\label{CramerResults}}

Using his probabilistic model of primes, Cram\'er 
obtained the following results for {\em random primes} $P_k$:

\begin{itemize}
\item The expected value of $\Pi(x)$ (the number of RPs $\le x$) is asymptotic to $\li x$: \ 
      $E\Pi(x) \sim \li x$ \ as $x\to\infty$.
\item With probability 1, we have: \ 
      $\displaystyle\limsup_{x\to\infty}{|\Pi(x)-\li x|\over\sqrt{2x\log\log x/\log x}}=1.$
\item With probability 1, we have: \ 
      $\displaystyle\limsup_{k\to\infty}{P_{k+1}-P_k\over\log^2{\negthinspace}P_k}=1.$
\end{itemize}
The latter result, {\em restated for true primes $p_k$}, constitutes the well-known {\em Cram\'er's conjecture}.
The conjecture appears likely to be true. Indeed, computations of Oliveira e Silva, Herzog \& Pardi (2014) 
have verified that
\begin{eqnarray*}
p_{k+1}-p_{k} &<& \log^2{\negthinspace}p_k \quad\mbox{for all primes }  p_k \notin \{2,3,7\} \mbox{ up to } 4\times10^{18},\\
p_{k}-p_{k-1} &<& \log^2{\negthinspace}p_k \quad\mbox{for all primes up to } 4\times10^{18}.
\end{eqnarray*}
At the same time, there exist large primes $p_k$ for which \ 
$0.8 < \displaystyle{p_{k+1}-p_k\over \log^2{\negthinspace}p_k} < 1$ \ (Sloane, 2014, A111943).


\section{Maximal gaps between Cram\'er's random primes}

\subsection{The exact distribution of maximal gaps}

We continue where Cram\'er left off. 
To obtain exact distributions of maximal gaps between Cram\'er's random primes, for now
we will restrict ourselves to finite sets of $n$ consecutive urns $U_1$, $U_2, \ldots, U_n$.
When our set of urns is small, we can compute the exact distributions of maximal gaps by hand, 
even without a computer.

For example, for $n=3$ we have only three urns: $U_1$, $U_2$, $U_3$.
Of these, only $U_3$ produces random results:
\begin{itemize}
\item a white ball (RP) with probability $\displaystyle{1\over \log3}\approx0.91$ \ {\em or}
\item a black ball (RC) with probability $\displaystyle1-{1\over \log3}\approx0.09$.
\end{itemize}
Thus, for the longest run of RCs, $R_3$, we have 
$R_3=0$ with probability 0.91, and $R_3=1$ with probability 0.09.
Consequently, for the maximal gap between RPs, we have $G_3=1$ with probability 0.91, and $G_3=2$ with probability 0.09.
(Recall that, by definition, $G_n = R_n+1$.)

One can visualize the exact distributions of maximal gaps in the form of histograms. 
With the help of a computer, we can find the exact distributions up to, say, $n=250$ urns.
Figure 1 shows the respective computer-generated distributions of maximal gaps
(cf.~Schilling, 1990, pp.\,197--204).

\begin{figure}[tbp] 
  \centering
  \includegraphics[bb=5 50 610 700,width=6.9in,height=6.9in,keepaspectratio]{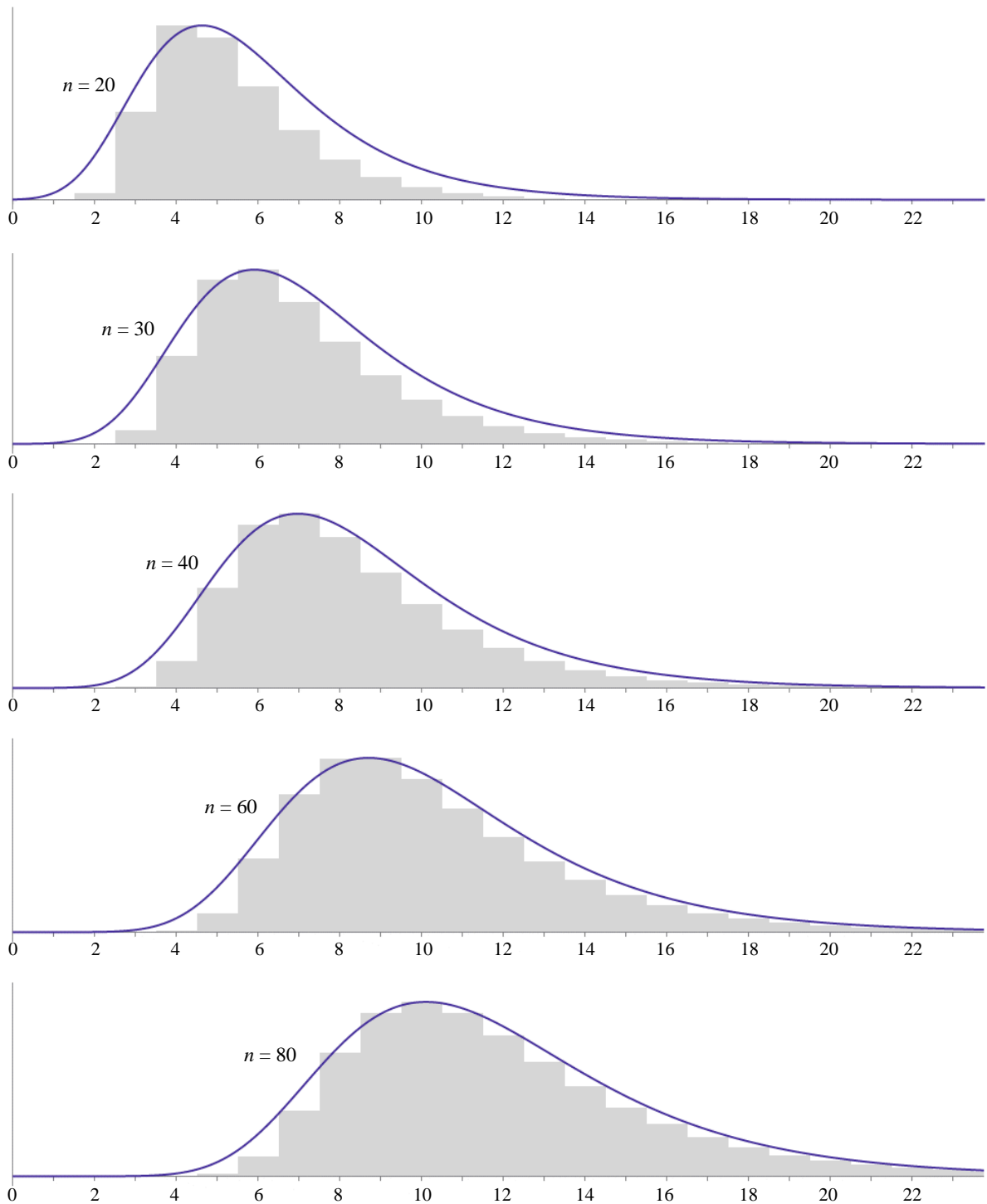}
  \caption{Histograms of the exact distributions of maximal prime gaps in Cramer's model with $n$ urns, computed for 
  $n=20,30,40,60,80$.
  The smooth curves represent Gumbel distributions (pdf) with the scale $a=n/\li n$ and mode $\mu=n\log(\li n)/\li n$.
  Interested readers can compute and plot the exact distributions of maximal gaps online at 
  {\small\tt http://www.javascripter.net/math/statistics/maximalprimegapsincramermodel.htm}.
  }
  \label{fig:fig1pdf}
\end{figure}

\medskip
\subsection{The limiting distribution of maximal gaps\label{TADOMG}}

In Figure 1, 
the exact distributions (histograms) of maximal gaps between RPs appear to approach the pdf curves 
of the Gumbel distribution with scale $\alpha_n = {n/\li n}$ and mode $\mu = \alpha_n \log\li n$.
We can restate this observation in a more precise form:

\textbf{Theorem 1}~\emph{
In Cram\'er's model with $n$ urns, the Gumbel distribution $\exp(-e^{-z})$ is 
the limiting distribution of maximal gaps $G_n$ between RPs: there exist $a_n>0$ and $b_n$ such that
$$
\lim_{n\to\infty} P(G_n \le x \equiv a_n z + b_n) = \exp(-e^{-z}), \quad\mbox{ where \ }
a_n \sim \alpha_n = {n\over\li n}, \quad
b_n \sim \alpha_n \log\li n.                     
$$
Equivalently, for $R_n$ (the longest runs of RCs) we have 
$\displaystyle\lim_{n\to\infty} P(R_n \le a_n z + b_n) = \exp(-e^{- z})$.  
}     

\smallskip\noindent
We will sketch two proofs of Theorem 1. The first proof will use the following lemmas.


\smallskip\noindent
{\bf Lemma of Common Median \ }
{\em
Suppose two Gumbel distributions have a common median and different scales $a\pm\varepsilon$, where 
$0<\varepsilon<{a\over2}$. Then the cdfs of these Gumbel distributions differ by no more than $\varepsilon\over a$. 
}

\smallskip\noindent
{\bf Lemma of Common Scale \ }
{\em
Suppose two Gumbel distributions have a common scale $a>0$ and medians $M\pm\delta$.
Then these Gumbel cdfs differ by no more than $\delta\over a$. 
}

\smallskip\noindent
Denote by $F_n(x)$ the cdf of the exact distribution of maximal gaps in Cram\'er's model with $n$ urns $U_1,\ldots,U_n$. 
Denote by $I_n$ ($n\ge10$) the largest interval of the $x$ axis such that $F_n(x)\in[(\log n)^{-1},1-(\log n)^{-1}]$ 
for all $x\in I_n$. (One can show that $I_n \subset [\log n,\log^2{\negthinspace}n]$.)

\smallskip\noindent
{\bf Squeeze Lemma \ }
{\em
Let $M_n$ be the median of $F_n(x)$. Then 
$F_n(x)$ is squeezed in the area bounded by the Gumbel distribution cdfs with medians $M=M_n\pm\delta$ 
and scales $a={n/\li n}\pm\varepsilon$:
$$
\min_{M\,=\,M_n\pm\delta \atop a\,=\,n/\li n\,\pm\,\varepsilon} 2^{-e^{-{x-M\vphantom{f}\over\vphantom{f}a}}}
~\le~ F_n(x) ~\le~ 
\max_{M\,=\,M_n\pm\delta \atop a\,=\,n/\li n\,\pm\,\varepsilon} 2^{-e^{-{x-M\vphantom{f}\over\vphantom{f}a}}}
\quad\mbox { for } \ n\ge10, \ \ x \in I_n,
$$
where one can take} $\varepsilon=3/2$, $\delta=1$.
(See Fig.\,2. 
The sequence $\{M_n\}$ is OEIS A235492, Sloane, 2014.)

\smallskip\noindent
Hereafter we often use these formulas expressing the Gumbel distribution cdf 
in terms of its scale $a$, mode $\mu$ and median $M$:
$$
\Gumbel(x;a,\mu) ~=~ e^{-e^{-{x-\mu\over\vphantom{f}a}}} ~=~ 2^{-e^{-{x-M\vphantom{f}\over\vphantom{f}a}}}, 
\quad\mbox{ where } \ M\,=\,\mu-a\log\log2\,\approx\,\mu+0.3665a.
$$

\subsubsection{First proof of Theorem 1}

Suppose the number of urns $n$ is large. We will examine two cases:
(a) $x \in I_n$; \ 
(b) $x \notin I_n$.

\medskip\noindent
{\bf Case (a) } First, we consider $F_n(x)$ for $x \in I_n$, i.\,e., $(\log n)^{-1}\le F_n(x)\le1-(\log n)^{-1}$.

\noindent
Let us estimate $M_n$, the median of $F_n(x)$. 
As before, let $\alpha_n=n/\li n$ and $\varepsilon=3/2$.
We can approximate Cram\'er's 
proportions of white balls in urns using two different procedures:

\noindent
(i) set the white balls' probability in all urns to $\displaystyle{1\over\alpha_n+\varepsilon}$, 
then increase the percentages of white balls in all urns \ {\em or, alternatively,}

\noindent
(ii) set the white balls' probability in all urns to $\displaystyle{1\over\alpha_n-\varepsilon}$, 
then reduce the percentages of white balls in all but a small subset of urns.

\smallskip\noindent
Note that the white balls probability ${1/a}$ can be approximated by the 
exponential distribution of gaps $\Exp(x;a)$.
Observe also that increasing the percentage of white balls pushes the median of maximal gaps to the left,
while reducing this percentage pushes the median to the right.
Therefore, $M_n$ must be somewhere between the medians of 
$(\Exp(x;\alpha_n-\varepsilon))^{\li n}$ and $(\Exp(x;\alpha_n+\varepsilon))^{\li n}$: 
$$
\mbox{ median } (\Exp(x;\alpha_n-\varepsilon))^{\li n} ~\lesssim~ M_n ~\lesssim~
\mbox{ median } (\Exp(x;\alpha_n+\varepsilon))^{\li n}.
$$
Since $\alpha_n\sim\log n$ while $\varepsilon=O(1)$,
the above lower and upper bounds are asymptotic to each other and 
to the median of $(\Exp(x;\alpha_n))^{\li n}$, so we must have
$$
M_n ~\sim~ \mbox{ median } (\Exp(x;\alpha_n))^{\li n} \quad\mbox{ as } n\to\infty.
$$ 
That is to say, as $n\to\infty$, the median $M_n$ of $F_n(x)$ must be asymptotic to 
the median of the limiting distribution of maxima of $\lfloor\li n\rfloor$ i.i.d.\,\,random variables 
with the exponential distribution $\Exp(x;\alpha_n)$.
But this limiting distribution is precisely the Gumbel distribution 
$\Gumbel(x;\alpha_n,\alpha_n\log\li n)$ 
(Gumbel, 1935; Gnedenko, 1943; Hall \& Wellner, 1979); therefore
$$
M_n ~\sim~ \mbox{median} \big(\Gumbel(x;\alpha_n,\alpha_n\log\li n)\big) 
    ~\approx~ \alpha_n\log\li n + 0.3665\alpha_n.
$$ 
On the other hand, it follows from Lemmas that $F_n(x)$ is squeezed in the 
$\Delta_n$-neighborhood of the Gumbel cdf $2^{-e^{-{x-M_n\over \alpha_n}}}$, where
$\varepsilon$ and $\delta$ are defined as in Squeeze Lemma,
$$
\Delta_n ~=~ {\varepsilon+\delta\over\alpha_n}, \qquad\mbox{ and }\qquad
\alpha_n ~=~  {n\over\li n}\sim \log n \to\infty \ \ \mbox{ as } n \to\infty.
$$
It is easy to see that $\Delta_n\to0$ as $n\to\infty$.

To sum it up: 
the median $M_n$ is asymptotic to the median of $\Gumbel(x;\alpha_n,\alpha_n\log\li n)$, 
while the limiting shape of $F_n(x)$ is dictated by the fact that $F_n(x)$ is in the $\Delta_n$-neighborhood 
of the Gumbel cdf\, $2^{-e^{-{x-M_n\over \alpha_n}}}$, with $\lim\limits_{n\to\infty}\Delta_n=0$.
This completes the proof for case (a).

\medskip\noindent
{\bf Case (b) } Now consider $x \notin I_n$: either $F_n(x)<(\log n)^{-1}$ or $F_n(x)>1-(\log n)^{-1}$.
Taking into account the monotonicity of cdfs, we can conclude from Lemmas that, for large $n$, 
\begin{eqnarray*}
  2^{-e^{-{x-M_n\over\vphantom{f}\alpha_n}}}&<&{\varepsilon+\delta\over\alpha_n}+(\log n)^{-1} 
  \qquad\mbox{ if }\mspace{37mu}F_n(x)~<~(\log n)^{-1},
\\
1-2^{-e^{-{x-M_n\over\vphantom{f}\alpha_n}}}&<&{\varepsilon+\delta\over\alpha_n}+(\log n)^{-1} 
  \qquad\mbox{ if }\ \ 1-F_n(x)~<~(\log n)^{-1},
\end{eqnarray*}
where $\alpha_n=\displaystyle{n\over\li n}\sim \log n$. 
Therefore, $\big|F_n(x)-2^{-e^{-{x-M_n\over\vphantom{f}\alpha_n}}}\big|=O\big((\log n)^{-1}\big)\to0$~~as $n\to\infty$. \hfill$\square$

\begin{figure}[tb] 
  \centering
  \includegraphics[bb=27 -2 583 239,width=6in,height=6in,keepaspectratio]{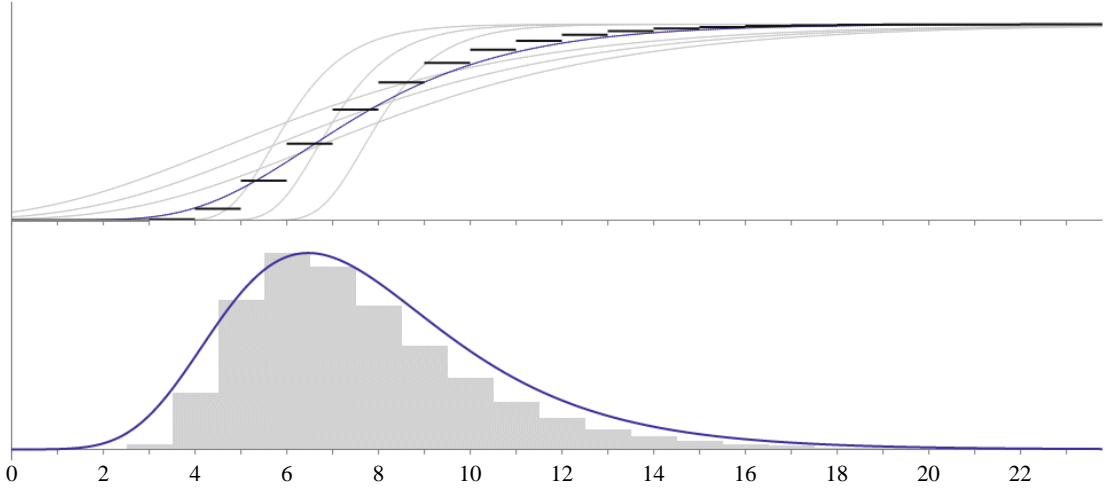}
  \caption{
   Bottom: the histogram of maximal prime gaps in Cram\'er's model with 35 urns.
   Top: the black ``staircase'' is the exact cdf of the distribution of maximal prime gaps in Cram\'er's model with $35$ urns.
   The dark curve shows the corresponding Gumbel distribution with scale $a=n/\li n$ and mode $\mu=n\log(\li n)/\li n$.
   Light gray curves illustrate the Gumbel cdfs with scales $a=n/\li n\pm\varepsilon$ and medians $M=M_n\pm1$ used in Squeeze Lemma; 
   $n=35$.
  }
  \label{fig:fig2cdf}
\end{figure}

\subsubsection{Second proof of Theorem 1}

We begin by proving the theorem for {\em longest runs $R_n$ of ``composites''}.
In Cram\'er's model, urns $U_n$ ($n\ge3$) produce white balls with a monotonically decreasing, slowly 
varying probability $1/\log n\to0$ as $n\to\infty$.
We observe that the more general Theorem A.1 (see {\it Appendix}) is applicable to our situation.
Theorem A.1 tells us that, as $n\to\infty$, the limiting distribution of longest runs {\em exists};
it is the Gumbel distribution with the scale and mode parameters determined by $E\Pi(n)$, 
the expected total number of RPs~$\le n$:
\begin{eqnarray*}
\mbox{scale } a_n &\sim& {n\over E\Pi(n)} ~=~ {n\over\li n + O(1)} ~\sim~ {n\over\li n} \\
\mbox{ mode } b_n &\sim& {n\log E\Pi(n)\over E\Pi(n)} ~=~  {n\over\li n + O(1)} \log(\li n +O(1)) ~\sim~ {n\over\li n} \log\li n.                 
\end{eqnarray*}
Here we have used the fact that, for $n\ge3$, the expected total number of RPs~$\le n$ is
$$
E\Pi(n) ~=~ 1+\sum_{k=3}^n {1\over\log k} ~=~ \li n + O(1) \quad\mbox{ as } n\to\infty. 
$$
Thus we can use the above scale and mode as rescaling parameters $a_n$ and $b_n$,
to obtain the standard Gumbel distribution $\exp(-e^{-z})$.
Note that $1\ll a_n\ll b_n$ as $n\to\infty$; so the distribution rescaling formula
$z = \displaystyle{x - b_n\over a_n}$ will produce approximately equal values of $z$, 
no matter whether we are rescaling the longest runs $R_n$ or maximal gaps $G_n\equiv R_n+1$.
\hfill$\square$

\subsection{Properties of the distribution of maximal gaps}

Let us look at the the properties of the distribution of maximal gaps.
We can readily see that the exact distribution of maximal gaps is {\em discrete} and {\em bounded},
while the limiting Gumbel distribution is {\em continuous}, {\em smooth}, and {\em unbounded}.
Below we discuss two properties that are common to the exact and asymptotic distributions of maximal gaps: 
{\em log-concavity} and {\em unimodailty}.

\subsubsection{Log-concavity}

A function $f(x)>0$ with a convex domain is {\em log-concave} if \ 
$\log f(\lambda x+(1-\lambda) y) \ge \lambda\log f(x) + (1-\lambda)\log f(y)$ \ 
for all $0\le\lambda\le1$ and for all $x, y$ in the domain of $f(x)$. 
%
A discrete sequence $\{s_k\}$ is log-concave if $s_k^2 \ge s_{k-1}s_{k+1}$ for each middle term $s_k$.

\noindent
{\bf Log-concavity of the limiting distribution of maximal gaps.} 
It is well known (and easy to check given the formulas for the distribution's pdf and cdf)
that the Gumbel distribution pdf and cdf are log-concave. Note that, in general, if $f(x)$ is a continuous 
distribution pdf and $F(x)$ is the corresponding cdf, then the following implications are true:
\begin{eqnarray*}
f(x) \mbox{ is log-concave} &\Rightarrow &   F(x) \mbox{ is log-concave}, \\
f(x) \mbox{ is log-concave} &\Rightarrow & 1-F(x) \mbox{ is log-concave}.
\end{eqnarray*}

\noindent
{\bf Log-concavity of exact distributions of maximal gaps.} 
A direct computational check shows that all the exact distribution functions we have computed in 
Section 4.1 are also log-concave:
$$
F_n(k)^2 ~\ge~ F_n(k-1)\,F_n(k+1), \qquad 1 < k < n.
$$
However, the log-concavity is not necessarily preserved
if we use an approximate (e.\,g.~Monte-Carlo) algorithm rather than the exact formulas for computing 
the finite distribution functions $F_n(x)$.

\subsubsection{Unimodality}

The Gumbel distribution $\Gumbel(x;a,\mu)$ is {\em unimodal} for any $\mu$ and any $a>0$:
it has a unique mode (most probable value), namely, $\mu$. 
What about the exact distribution of maximal gaps between RPs in Cram\'er's model with $n$ urns?
By Theorem 1, the Gumbel distribution is the limiting distribution of maximal gaps between RPs in Cram\'er's model;
therefore it is reasonable to expect that, for large $n$, 
in Cram\'er's model with $n$ urns the exact distribution of maximal gaps between RPs is also unimodal.
On the other hand, for moderate values of $n$ 
one can check by direct computation that each of the exact distributions of maximal gaps between RPs 
is unimodal. Here we will state without a formal proof the following  

{\bf Unimodality Lemma \ }
{\em In Cram\'er's model with $n$ urns, the exact distribution of maximal gaps between RPs is unimodal for each} $n>1$.

\section{Generalizations. Applications. Concluding remarks} 

We have thus shown that the Gumbel distribution is the limit law for extreme gaps between Cram\'er's random primes.
In particular, when the number of urns $n$ in Cram\'er's model is large, the distribution of
maximal prime gaps approaches the Gumbel distribution $\Gumbel(x;a_n,b_n)$ with these parameters:
$$
\mbox{scale } a_n ~\sim~ {n\over\li n}, \qquad
\mbox{ mode } b_n ~\sim~ {n\over\li n} \log\li n. 
$$
In view of this limit law for maximal gaps between Cram\'er's RPs $\le n$,
additional questions naturally arise: 

\begin{itemize}
  \item Is the Gumbel distribution, after proper rescaling, also an (almost sure) limit law for record prime gaps 
        observed in {\em a single infinite sequence} of Cram\'er's random primes? 
        Cf.~Resnick (1973) and Arnold {\em et al.} (1998, p.\,193).
  \item Is the Gumbel distribution, after proper rescaling, also the limit law for record gaps between {\em true primes}?
\end{itemize}

While the latter question appears particularly difficult,
it clearly suggests that our probabilistic result has a potential {\em application to number theory}.
Indeed, in the spirit of Cram\'er's program we quoted in {\em Introduction}, a number-theoretic result 
can be first obtained heuristically via probability considerations, and a rigorous number-theoretic proof could be found afterwards.
 
A straightforward generalization extends the Gumbel limit law to maximal gaps between 
{\em prime constellations}~--- dense clusters of consecutive primes with a repeatable pattern (see Table 1).
For this generalization, it is essential that we regard prime constellations near $x$
as random events occurring with a slowly varying probability decreasing to zero as $x\to\infty$.
Just as with extreme gaps between primes, here we do not have rigorous number-theoretic proofs~--- 
but we do have a convincing probabilistic argument explained below.

\vspace{0.3cm}
Table 1. Prime constellations. 
Conjectured probabilities to find a constellation starting at $x$ \ (cf.~Forbes, 2012) \vspace{-0.3cm}
\begin{center}
\begin{tabular}{lcc}\hline
Constellation type& Pattern of primes\vphantom{\large$1^1$}&Conjectured probability $C(\log x)^{-K}$ at $p\approx x$ \\
[0.2ex]\hline
Twin primes      & $\{p,\ p+2\}$ \vphantom{\large$1^1$}    & \hphantom{1}$1.32032\,(\log x)^{-2}$ \\
Prime triplets   & $\{p,\ p+2,\ p+6\}$                     & \hphantom{1}$2.85825\,(\log x)^{-3}$ \\
Prime triplets   & $\{p,\ p+4,\ p+6\}$                     & \hphantom{1}$2.85825\,(\log x)^{-3}$ \\
Prime quadruplets& $\{p,\ p+2,\ p+6,\ p+8\}$               & \hphantom{1}$4.15118\,(\log x)^{-4}$ \\
Prime quintuplets& $\{p,\ p+2,\ p+6,\ p+8,\ p+12\}$        &            $10.13179\,(\log x)^{-5}$ \\
Prime quintuplets& $\{p,\ p+4,\ p+6,\ p+10,\ p+12\}$       &            $10.13179\,(\log x)^{-5}$ \\
Prime sextuplets & $\{p,\ p+4,\ p+6,\ p+10,\ p+12,\ p+16\}$&            $17.29861\,(\log x)^{-6}$ \\\hline
\end{tabular}
\end{center}\vspace{-0.3cm}
\vspace{0.3cm}

The idea to model prime constellations by random events with predictable probabilities stems from the 
{\em K-tuple conjecture} of Hardy \& Littlewood (1922). 
Riesel (1994, pp.\,60--68) gives an accessible account of this conjecture.
Forbes (2012) published extensive numerical data on densest permissible prime constellations;
Table 1 is based on a small subset of Forbes' data.
Once we agree to treat prime constellations as random events with slowly varying probabilities given in Table 1,
the Gumbel limit law for maximal gaps between such events immediately follows from Theorem A.1 ({\em Appendix}).
Specifically, we can estimate the scale and mode in the limit law $\Gumbel(x;a_n,b_n)$ for a particular constellation 
by expressing them in terms of $E\Pi_c(n)$, the expected total count of constellations below $n$:
$$
\mbox{scale } a_n ~\sim~ {n\over E\Pi_c(n)}, \qquad
\mbox{ mode } b_n ~\sim~ {n\log E\Pi_c(n)\over E\Pi_c(n)}, \qquad  
E\Pi_c(n) ~\approx~ C \displaystyle\int_0^n{\negthinspace}{(\log x)^{-K} dx},           
$$
where the values of $C$ and $K$ are given in the last column of Table 1. 
Thus prime constellations can be simulated using a modified version of Cram\'er's model:
urn $U_n$ produces a white ball with probability $C(\log n)^{-K}$.
(Of course, it is also possible to simulate prime constellations by $K$ consecutive 
white balls in the classical Cram\'er model as described in Section 3.
However, the latter approach is inferior because it does not accommodate the known values 
of the {\em Hardy-Littlewood constants} $C$ in Table 1.)

The generalization of the Gumbel limit law to extreme gaps between prime sextuplets (Table 1, last line)
leads to a somewhat unexpected application: the {\tt riecoin} cryptocurrency (Riecoin.org, 2014).
The computational process of {\tt riecoin} mining consists in finding very large prime sextuplets
(over 300 bits in size). The Gumbel limit law for maximal gaps between prime sextuplets
therefore also describes the distribution of ``worst cases'' in terms of computational work in {\tt riecoin} mining.

To conclude our brief discussion of applications, let us quote the British mathematician
G.\,H.\,Hardy (1877--1947) who made spectacular breakthroughs in the study of prime numbers~--- 
and paradoxically cautioned against emphasizing the applications and practical usefulness of mathematics:
\begin{quotation}\noindent
A science is said to be useful if its development tends to accentuate the existing inequalities
in the distribution of wealth, or more directly promotes the destruction of human life.
The theory of prime numbers satisfies no such criteria. Those who pursue it will, if they are wise,
make no attempt to justify their interest in a subject so trivial and so remote, and will
console themselves with the thought that the greatest mathematicians of all ages have found in it 
a mysterious attraction impossible to resist.
(Hardy, 1915) 
\end{quotation}

\section{Appendix: Maximal intervals between random events occurring with slowly varying probability $\ell(t)\to0$}

In this appendix we give a general theorem whose particular case ($\ell(t)=1/\log t$\, for $t\ge3$)
has been used in Section~4.2.2.
We say that a function $\ell(t)>0$ is {\em slowly varying}
if it is defined for positive $t$, and 
$\displaystyle\lim\limits_{t\to\infty}{\ell(\lambda t)\over \ell(t)} = 1$
for any fixed $\lambda>0$ \ (de Haan \& Ferreira, 2006, p.\,362).

{\bf Theorem A.1}\
{\em
Consider biased coins with tails probability $\ell(k)$ at the $k$-th toss, $0<\ell(k)<1$,
where $\ell(t)$ is a smooth, slowly varying, monotonically decreasing function, and
$\lim\limits_{t\to\infty}\ell(t) = 0$.
%
Then, after a large number $n$ of tosses, the asymptotic distribution of
the longest runs of heads $R_n$ is the Gumbel distribution: there exist $a_n>0$ and $b_n$ such that
$$
\lim_{n\to\infty} P(R_n \le x \equiv a_n z + b_n) = \exp(-e^{-z}), \quad\mbox{ where \ }
a_n \sim {n\over E\Pi(n)}, \quad
b_n \sim {n\log E\Pi(n)\over E\Pi(n)}.
$$
Here $E\Pi(n)$ is the mathematical expectation of the total number of tails $\Pi(n)$
observed during the first $n$ tosses: }
$$
E\Pi(n) = \sum\limits_{k=1}^{n}\ell(k).
$$

It may be surprising that the asymptotic distribution does exist at all
for {\em discrete} events occurring with a slowly varying probability $\ell(t)\to0$.
In contrast, for a biased coin with a constant positive probability of tails,
the limiting distribution of the longest run of heads does {\em not} exist (Schilling, 1990, p.\,203).
%
In the lemmas below, we assume the conditions of Theorem~A.1
and take $\lambda$ to be arbitrarily large ($\lambda\gg1$).
The lemmas are easy to prove using the theory of regularly varying functions
(see, e.\,g., de Haan \& Ferreira, 2006, pp.\,361--367).

\textbf{Lemma 1~}
$\displaystyle\lim_{t\to\infty}\Pi(t)=\infty$ a.s., \ while \
$\displaystyle\lim_{t\to\infty}{E\Pi(t)\over t}=0$.

\textbf{Lemma 2~}
$\displaystyle\lim_{t\to\infty}E(\Pi(\lambda t)-\Pi(t))=\infty$.

\textbf{Lemma 3~} $\displaystyle\lim_{t\to\infty}{(E\Pi(\lambda t) / E\Pi(t))}=\lambda$.

\textbf{Lemma 4~} $\displaystyle\lim_{t\to\infty}{E(\Pi(\lambda t)-\Pi(t)) \over E\Pi(\lambda t)}= {\lambda-1\over\lambda}$.

\smallskip\noindent
{\em Proof of Theorem A.1.}
Let $k$ denote the ordinal number of a coin toss, so the $k$-th toss has the tails probability $\ell(k)$.
Take $\lambda\gg1$ and consider the sequence $S(t,\lambda t)$ of consecutive tosses with $t<k\le\lambda t$.
For toss sequences $S(t,\lambda t)$ with larger and larger $t$, we see that:

\noindent
(i) $\ell(k)$ becomes nearly constant:
$\ell(\lambda t) \lesssim \ell(k) \lesssim \ell(t)$ \ (by monotonicity + slow variation);

\noindent
(ii) the expected total number of tails is $E(\Pi(\lambda t)-\Pi(t))\to\infty$ as $t\to\infty$ (by Lemma 2).

\noindent
As $t$ grows larger, the outcome of the toss sequence $S(t,\lambda t)$ becomes indistinguishable from
an equally long sequence of tosses of a {\em constant-bias} coin whose tails probability $q$ is
$$
q = {E(\Pi(\lambda t)-\Pi(t))\over\lambda t - t} \approx {E\Pi(\lambda t)\over\lambda t} \approx {E\Pi(n)\over n},
\quad\mbox{ where we set } n=\lfloor\lambda t\rfloor.
$$
In the latter setup with a constant-bias coin, head runs are modeled by
i.i.d.~geometric random variables (Schilling, 1990). The expected total number $m$ of head runs is
$$ 
m ~=~ q (\lambda t - t) ~=~ E(\Pi(\lambda t)-\Pi(t)) ~\approx~ {\lambda-1\over\lambda} E\Pi(n),
\quad\mbox{ with }n=\lfloor\lambda t\rfloor.
\eqno{(1)}
$$ 
In turn, these $m$ geometric i.i.d.~variables are highly accurately approximated
by exponential i.i.d.~variables 
whose common cdf is $\Exp(\xi;a)=1-e^{-\xi/a}$ (Anderson, 1970).
The scale parameter $a$ in $\Exp(\xi;a)$ is given by
$$
a ~=~ {\lambda t - t \over E(\Pi(\lambda t)-\Pi(t))} ~\approx~ {n\over E\Pi(n)} \qquad\mbox{(here and above we have used Lemma 4)}.
$$
The largest value $L$ of $m$ exponential i.i.d.\,random variables, with the 
cdf $1-e^{-\xi/a}$, has the limiting Gumbel distribution (Gumbel, 1935; Gnedenko, 1943; Hall \& Wellner, 1979):
$$
\lim_{m\to\infty} P(L \le a z + a \log m) = \exp(-e^{-z}).
$$
However, it is well known that for the geometrically distributed runs of heads with a constant-bias coin
(i.\,e.\,\,the tails probability $q=\mbox{const}$) the longest runs $R$ do not have a limiting distribution
(Anderson, 1970; Schilling, 1990):
$$ 
P(R \le a z + a \log m) = \exp(-e^{-z}) + O(q),
\eqno{(2)}
$$ 
with the geometric-to-exponential approximation error $O(q)$ 
preventing the convergence of the exact longest run distributions to the Gumbel distribution.
Nevertheless, in our original setup with slowly varying bias of the coin,
instead of a constant $q$ we have
$$ 
q_n \sim a_n^{-1} \sim \ell(n) \sim {E\Pi(n)\over n} \to0 \quad\mbox{ as } n\to\infty.
\eqno{(3)}
$$ 
For the sequence $S(t,\lambda t)$ to contain the absolute longest run of heads $R_n$ up to $n=\lfloor\lambda t \rfloor$,
we must take $\lambda$ very large, so ${\lambda-1\over\lambda}\approx1$, and from Equation~(1) 
we have $\log m \approx \log E\Pi(n)$.
%
Now Equations (2) and (3) 
yield the limiting distribution of $R_n$ under the theorem conditions:
$$
\lim_{n\to\infty} P(R_n \le a_n z + b_n) = \exp(-e^{-z}),
\quad \mbox{ where } a_n \sim {n\over E\Pi(n)}, \quad b_n \sim {n\,\log E\Pi(n)\over E\Pi(n)}.
$$
Thus the convergence to the limiting Gumbel distribution is restored. The convergence rate
clearly depends on the actual choice of the slowly varying function $\ell(t)\to0$.
\hfill$\square$

\smallskip\noindent
{\em Remarks.}
Theorem A.1 is general enough to describe maximal gaps for ``primes'' and ``prime constellations''
in Cram\'er's model. However, comparing it to the results of Mladenovi\'c (1999),
one might look for a further generalization.
For example, what kind of a limiting distribution (if any) would we get if in Theorem A.1
we replace the {\em slowly varying} probability $\ell(t)\to0$ by some
{\em regularly varying} probability $q(t)\to0$, $q(t)\in \mbox{RV}_\alpha$ with, say, $-1\le\alpha\le0$?
Here $\alpha$ is the index of regular variation; for definitions of $\alpha$ and  $\mbox{RV}_\alpha$
see, e.\,g., de Haan \& Ferreira (2006, p.\,362).
Note that $\alpha=0$ corresponds to slowly varying functions.

\pagebreak
\textbf{Acknowledgements}

The author expresses his sincere gratitude to the anonymous referees for their time and 
attention to this manuscript and for constructive criticism.
Special thanks also to Marek Wolf, Luis Rodriguez and Carlos Rivera
whose ideas and conjectures served as starting points for this work
(cf.\ Wolf, 2011; Rodriguez \& Rivera, 2009).

\smallskip\noindent
\textbf{References}

\hangafter=1
\setlength{\hangindent}{2em}
Anderson, C.~W. (1970). Extreme value theory for a class of discrete distributions with applications to some stochastic processes. 
{\it J.\ Appl.\ Prob.} {\it 7}, 99--113. http://dx.doi.org/10.2307/3212152

\hangafter=1
\setlength{\hangindent}{2em}
Arnold, B.~C., Balakrishnan, N., \& Nagaraja, H.~N. (1998). {\it Records}. New York, NY: {Wiley. \hfill} \linebreak 
http://dx.doi.org/10.1002/9781118150412

\hangafter=1
\setlength{\hangindent}{2em}
Caldwell, C. (2010). \emph{The gaps between primes}. Section 3:~Table and graph of record gaps.
Retrieved from http://primes.utm.edu/notes/gaps.html\#table

\hangafter=1
\setlength{\hangindent}{2em}
Cram\'er, H. (1936). On the order of magnitude of the difference between consecutive prime numbers. 
{\it Acta Arith.} {\it 2}, 23--46. http://matwbn.icm.edu.pl/ksiazki/aa/aa2/aa212.pdf

\hangafter=1
\setlength{\hangindent}{2em}
Fisher, R.~A., \& Tippett, L.\,H.\,C. (1928). Limiting forms of the frequency distribution of
the largest and smallest member of a sample, {\it Math.\ Proc.\ Cambridge\ Philos.\ Soc.} {\it 24}, {180--190. \hfill} \linebreak
http://dx.doi.org/10.1017/S0305004100015681

\hangafter=1
\setlength{\hangindent}{2em}
Forbes, A.~D.~(2012). \emph{Prime $k$-tuplets}. Section 21: List of all possible patterns of prime $k$-tuplets. 
The Hardy-Littlewood constants pertaining to the distribution of prime $k$-tuplets. {Retrieved from: \hfill} \linebreak http://anthony.d.forbes.googlepages.com/ktuplets.htm

\hangafter=1
\setlength{\hangindent}{2em}
Gnedenko, B.~V. (1943). Sur la distribution limite du terme maximum d'une s\'erie al\'eatoire. {\it Ann. Math.} {\it 44}, 423--453.
English translation: On the limiting distribution of the maximum term in a random series, in
S.~Kotz \& N.~L.~Johnson (Eds.) (1993).
{\it Breakthroughs in Statistics}, {\it 1}.
{\it Foundations and Basic Theory} (pp.\,185--225). New York, NY: Springer. 
http://dx.doi.org/10.2307/1968974

\hangafter=1
\setlength{\hangindent}{2em}
Granville, A. (1995). Harald Cram\'er and the distribution of prime numbers. 
{\it Scand.\ Actuar.\ J.} {\it 1}, 12--28{. \hfill} \linebreak 
http://dx.doi.org/10.1080/03461238.1995.10413946

\hangafter=1
\setlength{\hangindent}{2em}
Gumbel, E.~J. (1935). Les valeurs extr\^emes des distributions statistiques. 
{\it Annales de l'institut Henri Poincar\'e}, {\it 5}, 115--158.
http://archive.numdam.org/article/AIHP\_1935\_\_5\_2\_115\_0.pdf

\hangafter=1
\setlength{\hangindent}{2em}
Gumbel, E.~J. (1958). {\it Statistics of Extremes}. New York, NY: Columbia University Press. Mineola, NY: Dover.

\hangafter=1
\setlength{\hangindent}{2em}
de Haan, L., \& Ferreira, A. (2006). {\it Extreme Value Theory}. New York, NY: Springer. 

\hangafter=1
\setlength{\hangindent}{2em}
Hall, W.~J., \& Wellner, J.~A. (1979). The rate of convergence in law for the maximum of an exponential sample. 
{\it Statist.\ Neerlandica}, {\it 33}, 151--154. 
http://dx.doi.org/10.1111/j.1467-9574.1979.tb00671.x

\hangafter=1
\setlength{\hangindent}{2em}
Hardy, G.~H. (1915). Prime numbers. In \emph{Reports on the State of Science}. 
The 85th Meeting of the British Association for the Advancement of Science, Manchester, 
1915 (pp.\,350--355).
Reprinted in P.~Borwein, S.~Choi, B.~Rooney, A.\,Weirathmueller (Eds.) (2008).
\emph{The Riemann Hypothesis. A Resource for the Afficionado and Virtuoso Alike}. 
CMS Books in Mathematics, {\it 27} (pp.\,301--306). New York, NY: Springer. 

\hangafter=1
\setlength{\hangindent}{2em}
Hardy, G.~H., \& Littlewood, J.~E.~(1922). Some problems of `Partitio Numerorum.' III.
On the expression of a number as a sum of primes. \emph{Acta Math. 44}, 1--70.
http://dx.doi.org/10.1007/BF02403921

\hangafter=1
\setlength{\hangindent}{2em}
Kourbatov, A. (2013). Maximal gaps between prime $k$-tuples: a statistical approach.
{\it J.~Integer Sequences}, {\it 16}, Article 13.5.2. \ 
https://cs.uwaterloo.ca/journals/JIS/vol16.html \ \ 
http://arxiv.org/abs/1301.2242

\hangafter=1
\setlength{\hangindent}{2em}
Mladenovi\'c, P. (1999). Limit theorems for the maximum terms of a sequence of random variables 
with marginal geometric distributions. {\it Extremes}, {\it 2:4}, 405--419.
http://dx.doi.org/10.1023/A:1009952232519

\hangafter=1
\setlength{\hangindent}{2em}
Nevzorov, V.~B. (2001). {\it Records:~Mathematical Theory}, AMS Translations Series {\it 194}, Providence, RI: AMS. 

\hangafter=1
\setlength{\hangindent}{2em}
Nicely, T.~R. (2013). {\it First occurrence prime gaps}. Retrieved from
http://www.trnicely.net/gaps/gaplist.html


\hangafter=1
\setlength{\hangindent}{2em}
Oliveira e Silva, T., Herzog, S., \& Pardi, S. (2014).
Empirical verification of the even Goldbach conjecture and computation of prime gaps up to $4\cdot10^{18}$, 
{\it Math.~Comp.} {\it 83}, {2033-2060. \hfill} \linebreak 
http://www.ams.org/journals/mcom/2014-83-288/S0025-5718-2013-02787-1/

\hangafter=1
\setlength{\hangindent}{2em}
Pintz, J. (2007). Cram\'er vs.~Cram\'er: On Cram\'er's probabilistic model of primes.
{\it Funct.\ Approx.\ Comment.\ Math.} {\it 37}, 361--376.
http://dx.doi.org/10.7169/facm/1229619660

\hangafter=1
\setlength{\hangindent}{2em}
Resnick, S.~I. (1973). Record values and maxima. {\it The Annals of Probability}, {\it 1}, {650--662. \hfill} \linebreak
http://dx.doi.org/10.1214/aop/1176996892

\hangafter=1
\setlength{\hangindent}{2em}
Riecoin.org (2014). Riecoin: a cryptocurrency using prime sextuplets. Retrieved from http://www.riecoin.org

\hangafter=1
\setlength{\hangindent}{2em}
Riesel, H. (1994). \emph{Prime Numbers and Computer Methods for Factorization} (2nd ed.). Boston, MA: Birkh\"auser.

\hangafter=1
\setlength{\hangindent}{2em}
Rodriguez, L., \& Rivera, C. (2009).
Conjecture 66. Gaps between consecutive twin prime pairs.
Retrieved from 
http://www.primepuzzles.net/conjectures/      

\hangafter=1
\setlength{\hangindent}{2em}
Schilling, M.~F. (1990). The longest run of heads. {\it College Math. J.} {\it 21}, 196--207.
http://dx.doi.org/10.2307/2686886

\hangafter=1
\setlength{\hangindent}{2em}
Shanks, D. (1964). On maximal gaps between successive primes. {\it Math.~Comp.} {\it 18}, {646--651. \hfill} \linebreak
http://dx.doi.org/10.1090/S0025-5718-1964-0167472-8

\hangafter=1
\setlength{\hangindent}{2em}
Sloane, N.~J.~A.~(Ed.) (2014). {\it The On-Line Encyclopedia of Integer Sequences}. Published electronically {at: \hfill} \linebreak  
http://oeis.org/A235492

\hangafter=1
\setlength{\hangindent}{2em}
Wolf,~M. (2011). {\it Some heuristics on the gaps between consecutive primes}, preprint. Retrieved from{: \hfill} \linebreak  
http://arxiv.org/abs/1102.0481



\vspace{0.5cm}
\textbf{Errata fixed}

The following minor errata of previous versions have been corrected in arXiv:1401.6959v3 (this version).

In Section 4.2.1, ``white-to-black balls ratio'' should be ``white balls' probability''.

In Section 4.3.1, the 9th line should be: \ \ 
$f(x) \mbox{ is log-concave} ~\Rightarrow~ 1-F(x) \mbox{ is log-concave}.$

In Appendix, after equation (3): \ \ ${\lambda-1\over\lambda}\to1$ should be ${\lambda-1\over\lambda}\approx 1$
(passage to the limit was not intended or needed).

\vspace{0.5cm}
\textbf{Copyrights}

Copyright for this article is retained by the author(s), with first publication rights granted to the journal.

This is an open-access article distributed under the terms and conditions of the Creative Commons Attribution license (http://creativecommons.org/licenses/by/3.0/).

\end{document}